\documentclass[titlepage,12pt]{article}

\usepackage{amsfonts}
\usepackage{latexsym}
\usepackage[dvips]{graphicx}
\usepackage{epsf}
\begin{document}

\def\cite#1{{\bf{[#1]}}}
\newcommand{\del}{\nabla}
\newcommand{\calR}{{\cal{R}}}
\newcommand{\p}{\partial}
\newcommand{\HH}{{\Bbb{H}}}
\newcommand{\delR}{\del^{\R}}
\renewcommand{\H}{\HH}
\newcommand{\calP}{{\cal{P}}}
\newcommand{\calL}{{\cal{L}}}
\newcommand{\syst}{{\rm{syst}}}
\newcommand{\U}{{\cal{U}}}
\newcommand{\V}{{\cal{V}}}
\newcommand{\calH}{{\cal{H}}}
\newcommand{\calS}{{\cal{S}}}
\newcommand{\Probb}{{\hbox{Prob}}}
\newcommand{\calT}{{\cal{T}}}
\newtheorem{Definition}{Definition}[section]
\newtheorem{Def}{Definition}[section]
\newtheorem{Conj}{Conjecture}[section]
\newtheorem{Th}{Theorem}[section]
\newenvironment{Thnmb}{\vskip 6pt \noindent
{\bf{Theorem  \nmb}} \begin{it}}{\end{it} \vskip 6pt}
\newenvironment{Lmnmb}{\vskip 6pt \noindent{\bf{Lemma \nmb}}
\begin{it}}{\end{it} \vskip 6pt}
\newenvironment{quesnmb}{\vskip 6pt \noindent 
{\bf{Question \nmb}} \begin{it}} {\end{it} \vskip 6pt}
\newenvironment{Ack}{\vskip 6pt \noindent{\bf
Acknowledgements:} \quad}{\vskip 6pt}
\newenvironment{Remark}{\vskip 6pt \noindent{\bf
Remark:} \quad}{\vskip 6pt}
\newtheorem{Ex}{Example}[section]
\newenvironment{Excont}{\vskip 6pt \noindent{\bf Example 2.1,
Continued:}\quad}{\vskip 6pt}
\newenvironment{Exconti}{\vskip 6pt \noindent{\bf Example 2.3,
Continued:} \quad}{\vskip 6pt}
\newtheorem{Cor}{Corollary}[section]   
\newtheorem{Lm}{Lemma}[section] 
\newcommand{\Pf}{\noindent
{\bf{Proof}}:\quad} 
\newcommand{\Claim}{\noindent
{{\bf{Claim}}}: \quad} 
\newcommand{\ques}{\noindent
{{\bf{Question}}}: \quad} 
\newcommand{\Rc}{{\hbox{Ricci}}}  

\newcommand{\vol}{{\hbox{ \rm vol}}}
\newcommand{\Lp}{\Delta}
\newcommand{\dv}{{\hbox{div}}}
\newcommand{\grad}{{\hbox{grad}}}
\newcommand{\Z}{{\Bbb{Z}}}
\newcommand{\F}{{\Bbb{F}}}
\newcommand{\calF}{{\cal{F}}}
\newcommand{\calO}{{\cal{O}}}

\newcommand{\Q}{{\Bbb{Q}}}
\newcommand{\Fpn}{{\F}_{p^n}}
\newcommand{\const}{{\hbox{\rm const}}}
\newcommand{\tr}{{\hbox{tr}}}
\newcommand{\Vol}{{\hbox{Vol}}}
\newcommand{\area}{{\hbox{area}}}
\newcommand{\Slam}{{S_{\lambda}}}
\newcommand{\dist}{{\hbox{dist}}}
\newcommand{\intt}{{\hbox{int}}}
\newcommand{\eps}{\varepsilon}

\newcommand{\inj}{{\hbox{\rm inj. rad.}}}
\newcommand{\diam}{{\hbox{\rm diam}}}
\newcommand{\length}{{\hbox{length}}}
\newcommand{\Der}{{\hbox{Der}}}
\newcommand{\IDer}{{\hbox{IDer}}}
\newcommand{\res}{{\hbox{res}}}

\newcommand{\qed}{\noindent{\hbox{\vrule height1.5ex width.5em}}\vskip 6pt}

\newcommand{\Ricc}{{\hbox{Ricc}}}
\newcommand{\Lap}{\Delta}
\newcommand{\lam}{\lambda}
\newcommand{\mod}{{\hbox{mod}}\ }
\newcommand{\C}{{\Bbb{C}}}

\newcommand{\kap}{\kappa}
\newcommand{\ddt}{{{\p}\over{\p t}}}
\newcommand{\curly}{\preceq}

\newcommand{\btbt}{\left( \begin{array}{cc}}
\newcommand{\etbt}{ \end{array}\right)}
\newcommand{\bthbth}{\left( \begin{array}{ccc}}
\newcommand{\ethbth}{ \end{array}\right)}
\newcommand{\neps}{\not \eps}
\def\nin{{\not\in}}
\newcommand{\bcol}{\left(\begin{array}{c}}
\newcommand{\ecol}{\end{array}\right)}

\def\Res{{\hbox{Res}}}
\newcommand{\ind}{{\hbox{ind}}}
\newcommand{\Gal}{{\hbox{Gal}}}
\newcommand{\daggers}{\ddagger}
\def\bibitem#1{\item{\cite{#1}}}
\newcommand{\bibo}{\bibitem}

\def\char{{\hbox{char}}}
\def\plus{\oplus}

\newcommand{\bthree}{\left( \begin{array}{ccc}}
\newcommand{\ethree}{\ethbth}
\newcommand{\injrad}{{\hbox{\rm inj. rad.}}}
\newcommand{\genus}{{\hbox{\rm genus}}}
\newcommand{\calD}{{\cal{D}}}
\newcommand{\tauD}{\tau_{\calD}}

\title{Isoscattering on Surfaces}

\author{Robert Brooks\thanks{Partially supported by the Israel Science
  Foundation and the Fund for the Promotion of Research at the
  Technion}\\ Department of Mathematics\\
  Technion-- Israel Institute of Technology\\ Haifa, Israel  
\and Orit Davidovich\\ Department of Mathematics\\ Technion-- Israel
  Institute of Technology\\ Haifa, Israel}

\date{ \ {January, 2002}}

\maketitle

In this paper, we give a number of examples of pairs of non-compact
surfaces $S_1$ and $S_2$ which are isoscattering, to be defined
below. Our basic construction is based on a version of Sunada's
Theorem \cite{Su}, which has been refined using the technique of
transplantation (\cite{Be}, \cite{Zel}) so as to be applicable to
isoscattering. See \cite{BGP} and \cite{BP} for this approach, which
 is reviewed
below.
 
Our aim here is to present a number of examples which are
exceptionally simple in one or more senses. Thus, the present paper
can be seen as an extension of \cite{BP}, where the aim was to
construct isoscattering surfaces with precisely one end. We will show:

\begin{Th} \label{genus} 
\begin{description}

\item{(a)} There exist surfaces $S_1$ and $S_2$ of genus 0 with eight
  ends which are isoscattering.

\item{(b)} There exist surfaces $S_1$ and $S_2$ of constant curvature
  $-1$ which are of genus 0 and have fifteen ends.

\item{(c)} There exist surfaces $S_1$ and $S_2$ of genus $1$ with five
  ends, or genus $2$ with three ends, which are isoscattering.

\item{(d)} There exist surfaces $S_1$ and $S_2$ of constant curvature
  $-1$ which are of genus 1 with thirteen ends, or genus 2 
with five ends, or genus 3 with three
  ends, which are isoscattering.

\item{(e)} There exist surfaces $S_1$ and $S_2$ of genus 3 with one
  end, or constant curvature of genus 4 with one end, which are
  isoscattering. 
\end{description}

\end{Th}

Part (e) is just a statement of the results of \cite{BP}, and is
recorded here for the sake of completeness. It will not be discussed
further in this paper.

The nature of the ends in Theorem \ref{genus} is not too
important. In the cases where the curvature is variable, they 
can be taken to be hyperbolic funnels or Euclidean
cones, or to be hyperbolic finite-area cusps. In the constant
curvature $-1$ cases, they can be taken either to be infinite-area
funnels or finite-area cusps.

Recall that a surface $S$ is called a {\em congruence surface} if $S=
\H^2/\Gamma$, where $\Gamma$ is contained in $PSL(2, \Z)$ and contains
a subgroup 
$$ \Gamma_k = \left\{ \btbt a&b \\ c & d \etbt \equiv  \pm \btbt 1 &0\\ 0&1
\etbt (\mod k) \right\}$$ 
for some $k$. In other words, the group $\Gamma$ is the inverse image
of a subgroup of $PSL(2, \Z/k)$ under the natural map 
$$ PSL(2, \Z) \to PSL(2, \Z/k).$$

We then have:

\begin{Th} \label{cong} There exist two congruence surfaces $S_1$ and
  $S_2$ which are isoscattering.

\end{Th}

Theorem \ref{cong} answers a question which was raised to us by Victor
Guillemin. The point here is that congruence surfaces have a
particularly rich structure of eigenvalues embedded in the continuous
spectrum. On the other hand, subgroups of $PSL(2, \Z/k)$ have a very
rigid structure \cite{Di}, and it is not {\it a priori} clear that the
finite group theory is rich enough to support the Sunada method.

A version of Theorem \ref{genus} was announced without proof in an
appendix to \cite{BJP}.

Theorems \ref{s4} and \ref{D} are certainly well-known to finite-group
theorists. We hope that the explicit treatment given here will be
useful to spectral geometers.

\begin{Ack} The first author would like to thank MIT for its warm
  hospitality and the Technion for its sabbatical support for the
  period in which this paper was written. He would also like to thank
  Victor Guillemin for his interest, and Peter Perry for his
  suggestion to pursue these questions in the context of the paper
  \cite{BJP}.

\end{Ack}

\section{Transplantation and Isoscattering}

We recall the approach to the Sunada Theorem given in
\cite{BGP}. Recall that a Sunada triple $(G, H_1, H_2)$ consists of a
finite group $G$ and two subgroups $H_1$ and $H_2$ of $G$ satisfying

\begin{equation}\label{dagger}
{\hbox{for all}}\ g \in G, \#([g] \cap H_1 ) = \#([g] \cap H_2),
\end{equation}
where $[g]$ denotes the conjugacy class of $g$ in $G$.

\begin{Th}[{\cite{Su}, \cite{BGP}}]\label{sun}
 Suppose that $M$ is a manifold and
  $\phi: \pi_1(M) \to G$ a surjective homomorphism.

Let $M^{H_1}$ and $M^{H_2}$ be the coverings of $M$ with fundamental
groups $\phi^{-1}(H_1)$ and $\phi^{-1}(H_2)$ respectively. 

Then there is a linear isomorphism $$\calT: C^{\infty}(M^{H_1}) \to
C^{\infty}(M^{H_2}) $$
such that $\calT$ and $\calT^{-1}$ commute with the Laplacian.

\end{Th}

We first remark that condition (\ref{dagger}) is equivalent to the
following: if we denote by $L^2(G/H_i)$ the $G$-module of functions on
the cosets $G/H_i$, then

\begin{equation}
L^2(G/H_1)\ {\hbox{is isomorphic to}}\ L^2(G/H_2)\ {\hbox{as}}\
G-{\hbox{modules}}.
\end{equation}

We may further rewrite this by noting that $L^2(G)$ has two
$G$-actions, on the left and on the right, so that we may write
$$ L^2(G/H_i) \equiv (L^2(G))^{H_i},$$
where the equivariance under $H_i$ is taken with respect to the left
$G$-action, and $G$-equivariance is taken with respect to the right
$G$-action. Equation (\ref{dagger}) is then equivalent to

\begin{equation}\label{dags}
(L^2(G))^{H_1}\ {\hbox{is $G$-isomorphic to}}\ (L^2(G))^{H_2}.
\end{equation}

We may further rewrite this equation as saying that there is a
$G$-equivariant map $T: L^2(G) \to L^2(G)$ which induces an
isomorphism $(L^2(G))^{H_1} \to (L^2(G))^{H_2}$. 

Now any $G$-map is determined by its value on the delta function,
which in turn can be described by a function 
$$c: G \to R,$$
so that the $G$-module map is given by
$$T(f)(x) = \sum_{g \in G} c(g) f(g\cdot x).$$

The requirement that the image of this map lies in $(L^2(G))^{H_2}$
can be expressed in terms of $c$ by the condition that
\begin{equation}\label{star}
c(gh)= c(g)\ \quad {\hbox{for}}\ h \in H_2.
\end{equation}
We may therefore express the condition (\ref{dagger}) as the existence
of a function $c$ on $G$ which satisfies (\ref{star}), and which
furthermore induces an isomorphism as in (\ref{dags}). 

Given such a function $c$, we may then write out the function $\calT$
as follows: let $M^{id}$ be the covering of $M$ whose fundamental
group is $\phi^{-1}(id)$. Then we may identify $C^{\infty}(M^{H_i})$
with $(C^{\infty}(M^{id}))^{H_i}$. The desired expression for $\calT$
is then given by

$$\calT(f)(x) = \sum_g c(g) f(g\cdot x).$$

We emphasize that all of this makes perfectly good sense for any
function $c$ satisfying (\ref{star}). The condition that it induces an
isomorphism is the crucial property we need.

Clearly, $\calT$ and its inverse take smooth functions to smooth
functions, and also commute with the Laplacian, since both statements
are true of the action by $g$ and taking linear combinations.

This establishes the theorem.

We now consider the case when the manifold $M$ is complete and
non-compact. We will discuss here the case where $M$ is hyperbolic
outside of a compact set, the case of Euclidean ends having been
discussed in \cite{BP}. 

We begin with a complete surface $M_0$, and
consider a conformal compactification of $M_0$, consisting of one
circle for each funnel and a point for each cusp. We also pick a {\em
  defining function} $\rho$ on $M_0$, that is, a function which is
positive on $M$ and vanishes
to first order on the boundary of $M_0$. 

If $0 < \lambda <1/4$, then we choose real $s$ so that 
$$\lambda = (s)(1-s).$$
Then, if $f \in C^{\infty}(\p M_0)$, there exist unique functions $u$ on
$M_0$ and $\calS_s(f)  \in C^{\infty}(\p M_0) $ such that

\begin{description}

\item{(i)} $ \Lap(u) = \lambda u.$

\item{(ii)} $u \sim (\rho)^s \calS_s(f) + (\rho)^{1-s}f + \calO(\rho)$ 
as $\rho \to 0$.

\end{description}

The operator $\calS_s$ is the {\em scattering operator} for $s$, and
continues for all $s$ to be a meromorphic operator.  Two surfaces
$M_0$ and $M_1$ will be {\em isoscattering} if they have poles of the
same multiplicity at the same values of $s$. 

We now pick $M$ as in Theorem \ref{sun}, and lift the defining
function $\rho$ on $M$ to defining functions on $M^{H_1}$, $M^{H_2}$,
and $M^{id}$. Note that we may identify $\calS_s$ on $M^{H_i}$ with
the operator $\calS_s$ on the $H_1$-invariant part of $M^{id}$. We
then have:

\begin{Th} The surfaces $M^{H_1}$ and $M^{H_2}$ are isoscattering.

\end{Th}

\Pf If we are given $f \in C^{\infty}(\p M^{id})$, then clearly 
$$ \calT(u) \sim (\rho)^s \calT(S_f(f)) + (\rho)^{1-s} \calT(f) +
\calO(\rho),$$
or, in other words,

$$ \calT(\calS_s(f)) = \calS_s(\calT(f)).$$

Thus, $\calT$ intertwines $\calS_s$ for all $s$, and hence $\calS_s$
on $M^{H_1}$ and $M^{H_2}$ have poles (with multiplicities) at the
same values.

This completes the
proof.

\section{The Group $PSL(3, \Z/2)$}

It is a rather remarkable fact that most of the examples of
isospectral surfaces (\cite{BT}) as well as all of the examples of
Theorem \ref{genus}, can be constructed from one Sunada triple. This
is the triple $(G, H_1, H_2)$ , where
$$G= PSL(3, \Z/2),$$
and 
$$ H_1 = \bthree {}* & {}* & {}* \\ 0 & {}* & {}* \\ 0 & {}* & {}* \ethree
\quad
H_2 = \bthree {}* & 0 & 0\\ {}* & {}* & {}* \\ {}* & {}* & {}* \ethree.$$

Note that the outer automorphism 
$$ A \to (A^{-1})^t$$
takes $H_1$ to $H_2$, and also takes elements of $H_1$ to conjugate
elements. This is enough to show that $(G, H_1, H_2)$ is a Sunada
triple.

In this section, we will present the necessary algebraic facts to
prove Theorem \ref{genus}. Many of these facts are proved easily by
noting the isomorphism
$$PSL(3, \Z/2) \cong PSL(2, \Z/7).$$
It is somewhat difficult to see the subgroups $H_1$ and $H_2$ in
$PSL(2, \Z/7)$. The outer automorphism which takes $H_1$ to $H_2$ is,
however, easy to describe. It is the automorphism
$$\btbt a& b\\ c& d \etbt \to \btbt a & -b\\ -c & d \etbt = 
\btbt -1 & 0 \\ 0 & 1 \etbt \btbt a & b \\ c & d \etbt \btbt -1 & 0 \\
0 & 1 \etbt.$$
The fact that this cannot be made an inner automorphism follows from
the fact that $-1$ is not a square $(\mod 7)$.

We now describe the conjugacy classes of $PSL(2, \Z/7)$:

\begin{Lm} Every element of $PSL(2, \Z/7)$ is of order $1, 2, 3, 4$,
  or $7$. 

\begin{description}

\item{(a)} The only element of order 1 is the identity.

\item{(b)} Every element of order $2$ is conjugate to $\btbt 0 & 1\\
  -1 & 0 \etbt$.

\item{(c)} Every element of order $3$ is conjugate to $\btbt 1 & 1\\
  -1 & 0 \etbt$.

\item{(d)} Every element of order $4$ is conjugate to $\btbt 2 & 1 \\
  1 & 1 \etbt$.

\item{(e)} Every element of order $7$ is conjugate to either $\btbt 1
  & 1\\ 0 & 1 \etbt$ or $\btbt 1 & -1 \\ 0 & 1 \etbt$.

\end{description}

\end{Lm}

Translating back into the group $PSL(3, \Z/2)$ gives

\begin{Lm} The elements of $PSL(3, \Z/2)$ satisfy:

\begin{description}\label{con}

\item{(a)} Every element of order $2$ is conjugate to $$\bthree 1 & 1 &
  0\\ 0&1&0 \\ 0&0&1 \ethree.$$

\item{(b)} Every element of order $3$ is conjugate to 
$$\bthree 0&1&0 \\ 0&0&1\\ 1&0&0 \ethree.$$

\item{(c)} Every element of order $4$ is conjugate to
$$\bthree 1&1&0 \\ 0&1&1 \\ 0&0&1 \ethree.$$

\item{(d)} Every element of order $7$ is conjugate to either
$$\bthree 1 & 1 & 1 \\ 1& 1&0 \\ 0 & 1 & 1 \ethree$$
or
$$\bthree 1&0&1 \\ 1&1&1\\ 1&1& 0 \ethree.$$
\end{description}

\end{Lm}

\Pf It suffices to check that each matrix has the order indicated. We
remark that a simple criterion for an element to be of order $7$ is
that adding 1 to the diagonal entries produces a non-singular matrix. 

To check that the two matrices in (d) above are not conjugate, we
observe that their characteristic polynomials are distinct.

Identifying $G/H_1$ as non-zero row vectors, we now may calculate the
action of an element of $G$ on $G/H_1$ as a permutation
representation. We will be interested in the cycle structure of this
representation, which clearly only depends on the conjugacy class of
the element. It follows from the above that the same calculation is
also valid for the permutation representation on $G/H_2$.

\begin{Th} Let $g \in PSL(3, \Z/2)$. Then the cycle structure of the
  permutation representation of $g$ on $G/H_1$ and $G/H_2$ is given
  by:

\begin{description}

\item{(a)} If $g$ is of order $2$, then $g$ acts as the product of two
  cycles of order 2 and three 1-cycles.

\item{(b)} If $g$ is of order $3$, then $g$ acts as two $3$-cycles and
  a $1$-cycle. 

\item{(c)} If $g$ is of order $4$, then $g$ acts as a $4$-cycle, a
  $2$-cycle, and a $1$-cycle.

\item{(d)} If $g$ is of order $7$, then $g$ acts as a $7$-cycle.

\end{description}
\end{Th}

The proof is just an evaluation in each case of the representatives in
Lemma \ref{con}.

\section{Proof of Theorem \ref{genus}}

In this section, we will prove Theorem \ref{genus}. Our method will be
to find an orbifold surface $M$ and a surjective homomorphism from
$\pi_1(M)$ to a Sunada triple $(G, H_1, H_2)$. We will take $G$ to be
$PSL(3, \Z/2) \cong PSL(2, \Z/7)$, and $H_1$ and $H_2$ the
corresponding subgroups. We then would like to study the corresponding
coverings $M^{H_1}$ and $M^{H_2}$.

Let $x$ be a singular point of $M$, and let $g_x$ be the element of
$\pi_1(M)$ corresponding to going one around $x$, which is
well-defined up to conjugacy. Then the points of $M^{H_1}$ (resp.
$M^{H_2}$) lying over $x$ are in 1-to1 correspondence with the cycle
decomposition of $g_x$ on $G/H_1$ (resp. $G/H_2$). If $g_x$ acts
freely on the cosets, then we may choose the orbifold singularity at
$x$ so that it smooths out to a regular (i.e. nonsingular) point in
$G/H_1$ (resp. $G/H_2)$. By Theorem \ref{con}, this will happen only
if $g_x$ is the identity or of order $7$, in which case there will be
precisely one point lying over $x$. 

We wish to calculate the genus of $M^{H_1}$ and $M^{H_2}$ as
topological surfaces (that is, forgetting the orbifold structure). Our
strategy will be the following: we remove all the singular points of
$M$ and the points lying over them, so that the covering is now a
regular (i.e. non-orbifold) covering. We then multiply the Euler
characteristic of $M$ with the points removed by $7$, the index
$[G:H_1] = [G:H_2]$. We then add back in the points lying over the
singular points. By Theorem \ref{con}, there will be five such points
if $g_x$ is of order 2, three such points if $g_x$ is of order 3 or
order 4, and
one if $g_x$ is of order $7$. We may now compute the genus of
$M^{H_1}$ (resp. $M^{H_2}$) from the Euler characteristic.

We now have to calculate the number of ends. We must throw out those
singular points lying over a singular points of oreder 2 (there are
five of these), order 3, or order 4 (there are three of these in both
cases). We need not throw out the point lying over a singular point of
order $7$.

We then have to worry about whether $M^{H_1}$ and $M^{H_2}$ are
distinct. As argued, for instance, in \cite{Su} or \cite{BGP}, that if we
choose a variable metric on $M$, then for a generic choice of such
metrics $M^{H_1}$ and $M^{H_2}$ will be non-isometric. If we want
constant curvature metrics, this will in general fail if $M$ is a
sphere with three singular points, as examples in \cite{BT} show, but
if $M$ is a  sphere with $n$ singular points, $n \ge 4$, or a surface
of higher genus with an arbirary number of singular points, then
choosing a generic conformal structure on $M$ and generically placed
points will produce $M^{H_1}$ and $M^{H_2}$ distinct.

Let us first take the case where $M$ is a sphere with three singular
points. 

Choosing these singular points to be of order $2, 3$, and $7$, we must
find matrices $A, B,$ and $C$ in $PSL(2, \Z/7)$ such that:

\begin{description} 

\item{(i)} $A$ is of order 2, $B$ is of order $3$, and $C$ is of
  order $7$. 

\item{(ii)} $ABC= {\hbox{id}}$.

\item{(iii)} $A$, $B$, and $C$ generate $PSL(2, \Z/7)$.

\end{description}

A simple choice is
$$ A = \btbt 0&1\\ -1 & 0 \etbt \quad  B= \btbt 1 & 1\\ -1 & 0 \etbt$$
$$C= \btbt 1&0 \\ -1& 1 \etbt.$$

The computation of the genus of $M^{H_1}$ (resp. $M^{H_2}$) proceeds
as follows: the thrice-punctured sphere has Euler characteristic $\chi
= -1$. Hence $M^{H_i}$ without the singular points has $\chi =
-7$. Putting in the five singular points lying over the singular point
of order $2$, the three singular points lying over the point of order
$3$, and the point lying over the point of order $7$ adds $5 + 3 +
1=9$ to this, yielding an Euler characteristic of $2$. Hence $M^{H_1}$
and $M^{H_2}$ are of genus $0$. We must make ends out of the points
lying over the singular points of orders $2$ and $3$, to give a total
of eight ends. 

This establishes Theorem \ref{genus} (a).

If we had used two singular points of order $7$ and one singular point
of order $2$ (resp. 3 or 4), we would obtain for $M^{H_1}$ and $M^{H_2}$
surfaces of genus 1 (resp. 2) with five (resp. 3) ends, provided we
can find the corresponding generators. But 

$$B = \btbt 1 & k\\ 0 & 1 \etbt \quad {\hbox{and}}\ C= \btbt 1 & 0\\ l &
1 \etbt$$
generate $PSL(2,\Z/7)$ for any choice of $k,l$ prime to $7$, and their
product has trace $2+ kl$. Hence, for appropriate choice of $k$ and
$l$, we may find $A$ of order $2, 3,$ or $4$. 

This establishes (c).

We now investigate what happens when we choose $M$ to be a sphere with
four singular points.

Choosing the singular points to be of order $2, 2, 2$, and $7$
respectively, we calculate the Euler characteristic of $M^{H_i}$ as 
$$\chi = 7(-2) + 3\cdot 5 + 1 = 2,$$
so again the $M^{H_i}$ have genus $0$, now with fifteen ends, provided
we can find matrices $A, B, C,$ and $D$ of these orders which generate
$PSL(2, \Z/7)$ and whose product is 1. To do this, we first observe
that we may find two matrices $B'$ and $C'$ of order $2$ such that
their product is of order $3$. One choice is
$$B'= \btbt 0 & 1 \\ -1 & 0 \etbt , \quad C'= \btbt 0 & 2 \\ 3 & 0
\etbt.$$

We then conjugate $B'$ and $C'$ so that their product is $\btbt 1 & 1
\\ -1 & 0 \etbt$. We may then choose $A$ and $D$ to be $\btbt 0 & 1\\
-1 & 0 \etbt$ and $\btbt 1 & 0 \\ -1 & 1 \etbt$ as above.

This establishes (b).

To establish the first part of (d), we search for matrices of order
$2,2,3$, and $7$ whose product is 1. Choosing 
$$C = \btbt 1 & 1 \\ -1 & 0 \etbt, \quad D= \btbt 1 & 0\\ k & 1 \etbt,$$
we have that $C$ and $D$ generate $PSL(2, \Z/7)$, and for an
appropriate choice of $k$ the product is of order $3$. We may then
choose $A$ and $B$ as above to be two matrices of order two whose
product of order $3$ is the inverse of this matrix.

To establish the second and third parts of (d), we proceed
differently. The base surface $M$ will be of genus 1 with one singular
point. If the singular point is of order $2$, the resulting $M^{H_i}$
will be of genus $2$, with five ends. If the singular point is of
order $3$ or $4$, the resulting surface is of genus $3$, with three
ends. 

We therefore seek matrices $A$ and $B$ which generate $PSL(2, \Z/7)$,
such that their commutator is of order $2$ (resp. $3$ or $4$).

Choosing 
$$A= \btbt 0 & 1\\ -1 & 0 \etbt, \quad B = \btbt 1 & 0\\ k& 1 \etbt,$$
we get
$$\left[A,B\right]= \btbt 1 + k^2 & -k\\ -k & 1\etbt ,$$
and $A$ and $B$ generate $PSL(2,\Z/7)$. Choosing, for instance, $k=2$
gives a commutator of order $3$.

Choosing 
$$A = \btbt 0 & 1\\ -1 & 0 \etbt \quad B= \btbt 4& 1\\ 0 & 2 \etbt$$ 
gives
$$\left[A, B\right] = \btbt 3 & 2\\ 2&4 \etbt,$$
which is of order $2$.

We see no elegant way of seeing that $A$ and $B$ generate $PSL(2,
\Z/7)$, but
$$A(BAB^2)(\left[A,B\right])= \btbt 1& 0\\ 2 & 1 \etbt,$$
and this matrix and $A$ generate.

This concludes the proof of the second and third parts of (d), and
hence Theorem \ref{genus}.

\section{Proof of Theorem \ref{cong}}

In this section, we construct congruence surfaces which are
isoscattering.

There are several difficulties in this setting which are not present
in the general setting. First of all, congruence surfaces are
constructed out of subgroups of the finite groups $PSL(2, \Z/k)$, and
such groups are rather special. The subgroups of $PSL(2, \Z/p)$ have
been classified, and are given in Dickson's List \cite{Di}. It is not
{\em a priori} evident, for instance, that $PSL(2,\Z/p)$ contains
Sunada triples for general $p$. Of course, the case of $PSL(3, \Z/2)
\cong PSL(2, \Z/7)$ occurs as a very special example, but we will
need a richer collection of examples. 

Secondly, given such a Sunada triple $(G, H_1,H_2)$, we do not have
the freedom of choosing a homomorphism $\pi_1(M) \to G$ as we did
previously. It is given to us canonically. 

Finally, we must worry about ``extra isometries,'' since we do not
have the freedom to change parameters to guarantee that $M^{H_1}$ and
$M^{H_2}$ will be distinct.

We begin our discussion with the group $G= PSL(2, \Z/7)$, and $H_1$ and
$H_2$ the two subgroups as above. Taking
$$\Gamma= PSL(2, \Z),$$ 
and considering the natural projection $\phi: \Gamma \to G$, we first
notice that the $\phi^{-1}(H_i)$ contain torsion elements, so that the
$\H^2/\phi^{-1}(H_i)$ are singular surfaces.

To remedy this problem, and also to introduce a technique we will use
later, we note that for $k=k_1 k_2$, with $k_1$ and $k_2$ relatively prime,
we have
$$PSL(2, \Z/k)= P( SL(2, \Z/k_1) \times SL(2, \Z/k_2)).$$
Choosing $k=14$, we see that
$$PSL(2, \Z/14)= PSL(2,\Z/2)\times PSL(2,\Z/7),$$
noting that the ``$P$'' in $PSL(2, \Z/2)$ is trivial.

Furthermore, the kernel $\Gamma_2$ of $PSL(2, \Z/2)$ satisfies that
$\H^2/\Gamma_2$ is a (non-singular) thrice-punctured sphere. Hence, we
resolve the issue of singularities by restricting to $\Gamma_2$.

Now let $S^{i}= \H^2/\Gamma^{H_i},$ where
$$\Gamma^{H_i} = \Gamma_2 \cap \phi^{-1}(H_i).$$
Then, as before, $S^{1}$ and $S^{2}$ are isoscattering. They are,
however, also isometric. This can be seen by noting that $H_1$ and
$H_2$ are conjugate under the automorphism
$$\tau: \btbt a& b\\ c& d \etbt \to \btbt a & -b \\ -c & d \etbt.$$
Furthermore, this $\tau$ induces an orientation-reversing isometry of
$H^2/\Gamma$ which is reflection in the line ${\hbox{Re}}(z)=0$ in the usual
fundamental domain for $\H^2/\Gamma$, and therefore lifts to an
orientation-reversing isometry of $S^1$ to $S^2$. See the genus 3
example of \cite{BT}, where a similar problem occurs.

We will handle this problem in the following way: let us assume for
some $p$ different from $2$ or $7$, there is a subgroup $K$ of $SL(2,
\Z)$ such that $K$ and $\tau(K)$ are not conjugate in $SL(2,
\Z/p)$. We may now choose our subgroups

$$\begin{array}{ll}
G &= P({\hbox{id}}\times SL(2, \Z/7) \times SL(2, \Z/p))\\
H_1 &= P({\hbox{id}}\times H_1 \times K)\\
H_2 &= P({\hbox{id}}\times H_2 \times K),
\end{array}$$
and let $\widetilde{S}^1$ and $\widetilde{S}^2$ be the coverings of
$\H^2/PSL(2,\Z)$ corresponding to $\phi^{-1}(H_1)$ and
$\phi^{-1}(H_2)$. 

$\widetilde{S}^1$ and $\widetilde{S}^2$ are isoscattering, but we want to show
that they are not isometric. Any such isometry between them must be
given by conjugation by some matrix $C$, by the involution $\tau$, or
by a composition of the two.

The first possibility cannot obtain, because restricting to $PSL(2,
\Z/7)$, $C$ will give a conjugacy of $H_1$ to $H_2$. The second and
third possibilities also cannot hold, since if such a matrix $C'$
exists, restricting to the $PSL(2, \Z/p)$ factor, it will give a
conjugacy from $K$ to $\tau(K)$.

Thus, Theorem \ref{cong} will follow once we find such a $p$ and $K$.

We now examine Dickson's list \cite{Di} for likely subgroups $K$ of
$PSL(2, \Z/p)$ for which $K$ is not conjugate to $\tau(K)$. We remark
that $\tau$ is given by the outer automorphism
$$\btbt a&b\\ c&d \etbt \to \btbt -1 & 0 \\ 0 &1 \etbt \btbt a&b\\ c&d
\etbt \btbt -1 & 0 \\ 0 & 1 \etbt,$$
and for $p \equiv 1\ (\mod 4),$ we have that
 $-1$ is a square root $(\mod p)$, so
$\tau$ is actually an inner automorphism for such $p$.

To understand our choice of $K$, we observe that the subgroups $H_1$ and
$H_2$ are isomorphic to the symmetric group $S(4)$ on four elements 
(a more detailed discussion will be given below). We will show:

\begin{Th} \label{s4} For $ p \equiv 7\ (\mod 8)$, there exist
  subgroups $K$ of $PSL(2, \Z/p)$ isomorphic to $S(4)$, such that $K$
  is not conjugate to $\tau(K)$ in $PSL(2, \Z/p)$.
\end{Th}

\Pf We first discuss the restriction $p \equiv 7\ (\mod 8)$.

Indeed, $S(4)$ contains the cyclic subgroup $\Z/4$. In order for
$PSL(2, \Z/p)$ to contain a cyclic subgroup of order $4$, we must have
$p \equiv \pm 1 (\mod 8)$. The plus sign is ruled out by the condition
$p \not\equiv 1(\mod 4)$.

To construct a subgroup $K$ isomorphic to $S(4)$ in $PSL(2, \Z/p)$, we
first note that $S(4)$ contains the dihedral group
$$\{ A, D: A^4 =1, \quad  D^2=1, \quad DAD= A^{-1} \},$$
given by $A= (1,2,3,4)$, $D=(1,2)(3,4)$.

We now seek such a subgroup in $PSL(2, \Z/p)$. A convenient choice for
$A$ is
$$A = \btbt \alpha & \alpha\\ -\alpha & \alpha \etbt,\quad  \alpha^2
\equiv 1/2\ (\mod p).$$

Note that $2$ is a square $(\mod p)$ by the condition that $ p \equiv
7 (\mod 8)$ and quadratic reciprocity. Note also that
$$A^2= C_1 = \btbt 0&1\\- 1&0 \etbt.$$
We now seek a matrix $D$ such that $A$ and $D$ generate a dihedral
group of order 8. That means that

\begin{equation} \label{dihed}
D^2 = \btbt 1 &0 \\ 0&1 \etbt, \quad DAD= A^{-1}.
\end{equation}
The second relation implies that $D$ commutes with $C_1$. We observe
that any matrix $Z$ commuting with $C_1$ must either be of the form
\begin{equation}\label{ccomm1}
Z = \btbt x&y \\-y & x \etbt, \quad x^2 + y^2 =1
\end{equation}
or
\begin{equation}\label{ccomm2}
Z= \btbt \beta & \gamma\\ \gamma & -\beta \etbt, \quad \beta^2 +
\gamma^2 = -1.
\end{equation}
Since any matrix of the form (\ref{ccomm1}) commutes with $A$, $D$
must be of the form (\ref{ccomm2}). Furthermore, any matrix of the
form (\ref{ccomm2}) will satisfy (\ref{dihed}).

For later reference, we remark that neither $\beta$ nor $\gamma$ can
be zero, since $-1$ is not a square $(\mod p)$, and we may choose
$\beta + \gamma -1 \not \equiv 0 (\mod p)$, by changing the sign of
$\gamma$ if necessary.

We now set
$$C_2 = C_1 D = \btbt -\gamma & \beta \\ \beta & \gamma \etbt.$$
We may embed this dihedral group in $S(4)$ by setting
$$\begin{array}{ll}
A &\to (1,2,3,4)\\
C_1 &\to(1,3)(2,4)\\
D&\to (1,2)(3,4)\\
C_2 &\to (1,4)(2,3).
\end{array}$$

We now seek an element $E$ of $PSL(2, \Z/p)$ corresponding to the
element $(1,2,3)$. Thus, $E$ must satisfy the conditions

$$\begin{array}{ll}
E^3 &=1\\
EC_1E^{-1} &= D\\
EDE^{-1} &= C_2\\
EC_2E^{-1}&= C_1\\
AEA^{-1} &= E^{-1}C_1.
\end{array}$$
\newcommand{\Fix}{{\hbox{Fix}}}

To do this, let us tentatively set
$$\widetilde{E} = \btbt a&b\\ c&d \etbt.$$
We will want $\widetilde{E}$ to send the fixed points $\Fix(C_1)$ of
$C_1$ (viewed as a linear fractional tranformation) to $\Fix(D)$, the
fixed points $\Fix(D)$ of $D$ to $\Fix(C_2)$, and $\Fix(C_2)$ to $\Fix(C_1)$.

We calculate:

$$\begin{array}{ll}
\Fix(C_1) &= \pm i,\\
\Fix(D) &= \frac{\beta \pm i}{\gamma}\\
\Fix(C_2) &= \frac{ -\gamma \pm i}{\beta},
\end{array}$$
where we have denoted by $i$ a square root of $-1$ in the field
$\F_{p^2}$.

Choosing the plus sign in each case, we seek $\widetilde{E}$
satisfying

$$
\begin{array}{ll}
\widetilde{E}(i) = \frac{\beta + i}{\gamma}\\
\widetilde{E}( \frac{\beta + i}{\gamma}) = \frac{-\gamma + i}{\beta}\\
\widetilde{E}(\frac{-\gamma+ i}{\beta}) = i.
\end{array}$$

After some tedious linear algebra, we find

$$ \widetilde{E} = \btbt \frac{\gamma + \beta^2}{\beta + \gamma - 1}
& \frac{-\beta + \beta \gamma - \gamma}{\beta + \gamma -1}\\
\frac{ \beta \gamma - 1}{\beta + \gamma -1} & \frac{\beta +
  \gamma^2}{\beta + \gamma -1} \etbt.$$

Note that with these choices of $a,b,c,$ and $d$, we have
$$ a+d =1, \quad c-b=1, \quad ad-bc=1.$$

The first relation assures that $\widetilde{E}$ is of order $3$, and
it is easily checked that the conjugacies of $C_1$, $C_2$, and $D$ by
$\widetilde{E}$ are as desired.

It remains to check the condition
$$A\widetilde{E}A^{-1} = \widetilde{E}^{-1}C_1.$$
This unappetizing calculation can be carried out by writing

$$A\widetilde{E} A^{-1} \widetilde{E} = \btbt \alpha^2 & 0\\ 0&
\alpha^2 \etbt
\btbt 1& 1\\ -1 &1 \etbt \btbt a&b\\ c&d \etbt \btbt 1&-1\\ 1&1 \etbt
\btbt a&b\\ c&d \etbt,$$
using $a +d = 1$ and $c-b =1$ to eliminate $c$ and $d$, and $ad-bc =1$
to replace the quadratic terms $a^2 + b^2$ by $a-b-1$.  We find that
$$A\widetilde{E}A^{-1}\widetilde{E} =C_1$$
as desired.

Setting $E= \widetilde{E}$,  we now have constructed a subgroup $G$ of
$PSL(2, \Z/p)$ isomorphic to $S(4)$.

We now must show that $\tau(G)$ is not conjugate to $G$ in
$PSL(2,\Z/p)$.

But if $\psi: G \to G$ is any isomorphism, we may assume, by replacing
$\psi$ by a conjugate, that
$$\psi(A) = A.$$

$\psi(C_1)$ must then be $C_1$, and after conjugating by $A$ if
necessary, we must have
$$\psi(D)= D, \quad \psi(C_2)= C_2.$$
It then follows that
$$\psi(E) =E.$$
We now show that under these assumptions, we cannot have
$$ \psi(X) = Z\tau(X) Z^{-1}, \quad X \in G.$$
Noting that $\tau(C_1)= C_1$, this gives us
$$ C_1 = Z C_1 Z^{-1},$$
so that $Z$ must be of the form (\ref{ccomm1})
or(\ref{ccomm2}). 
Noting that $\tau(A)= A^{-1}$, this gives
$$A = Z A^{-1}Z^{-1},$$
so that $Z$ must be of the form 
$$Z= \btbt x & y\\ y & -x \etbt, \quad x^2 + y^2 = -1.$$

We now consider the equation
$$ D= Z \tau(D) Z^{-1},$$
which we write out as
$$ \pm \btbt \beta & \gamma \\ \gamma& -\beta \etbt = \btbt x& y\\
y&-x \etbt \btbt -\beta & \gamma\\ \gamma & \beta \etbt \btbt x&y\\
y&-x \etbt.$$

The term on the right is computed to be
$$\btbt -x^2\beta + 2xy\gamma + y^2 \beta & -2xy\beta + y^2 \gamma -
x^2\gamma\\ -2xy\beta -x^2\gamma + y^2 \gamma & -y^2\beta -2xy\gamma +
x^2 \beta \etbt.$$

Taking first the plus sign, we get from the upper-left entry the equation
$$ -x^2 \beta + 2xy\gamma + y^2\beta = \beta = \beta(-x^2 -y^2),$$
or
$$2xy\gamma + 2y^2\beta=0.$$

Since $y \neq 0$, this gives
$$ x\gamma + y\beta =0.$$
Similarly, the lower-left entry gives
$$ -x\beta + y\gamma=0.$$
Solving these equations gives
$$\frac{x}{\gamma}(\gamma^2 + \beta^2)=0.$$
But this implies $x=0$, a contradiction.

Now we take the minus sign. We get the equations
$$ -x^2 \beta + 2xy\gamma + y^2\beta = -\beta = \beta(x^2 + y^2)$$
and
$$-x^2 \gamma -2xy\beta + y^2 \gamma = \gamma(x^2 + y^2).$$
These become the two equations
$$ y\gamma = x\beta$$
and
$$-y\beta = x\gamma,$$
which again yields a contradiction. 

This contradiction proves Theorem \ref{s4}, and hence also Theorem
\ref{cong}. 

Another approach to Theorem \ref{cong} may be based on $p \equiv 1 \
(\mod 4)$. For $\calD \in \Z/p$, let 
$$\tauD \btbt a&b\\ c&d \etbt = \btbt \calD & 0\\ 0& 1 \etbt \btbt
a&b\\ c&d \etbt \btbt \calD^{-1} & 0\\ 0&1 \etbt = \btbt a& \calD b\\
\frac{1}{\calD} c & d \etbt.$$

When $\calD$ is a square $(\mod p)$, $\tauD$ is inner. We will show:

\begin{Th}\label{D} Let $p \equiv 1\ (\mod 4)$, and let $\calD$ be a
  non-square $(\mod p)$. 

Then there exists a subgroup $K$ of $PSL(2,\Z/p)$ isomorphic to
$S(4)$, such that $K$ and $\tauD(K)$ are not conjugate in
$PSL(2,\Z/p)$.

\end{Th}

Given Theorem \ref{D}, we set
$$\Gamma^1 = \Gamma_2 \cap \phi^{-1}(K)$$
and
$$\Gamma^2= \Gamma_2 \cap \phi^{-1}(\tauD(K)).$$
Then $S^1 = \H^2/\Gamma^1$ and $S^2=\H^2/\Gamma^2$ are isoscattering. The
orientation-reversing isometry $\tau$ of $\H^2/PSL(2,\Z)$ lifts to an
isometry of $S^i$ to itself, $i=1, 2$, since $-1$ is a square $(\mod
p)$, so Theorem \ref{D} suffices to show that they are not isometric.

To prove Theorem \ref{D}, we set
$$A = \btbt \alpha(1+i) &0\\ 0 & \alpha(1-i) \etbt \quad C_1 = \btbt
i&0\\ 0&-i \etbt, $$
and
$$D= \btbt 0&1\\ -1 & 0 \etbt,$$
where $i$ is a square root of $-1$ $(\mod p)$, and $\alpha^2 =
\frac{1}{2}\ (\mod p)$. We may then solve for $E$ as above, to find
$$ E= \btbt \frac{1-i}{2} & \frac{1-i}{2}\\ -\frac{(1+i)}{2} &
\frac{1+i}{2} \etbt.$$

Then 
$$\tauD(A)= A, \quad \tauD(C_1)= C_1,$$
and
$\tauD(D) = \btbt 0 & \calD\\ - \frac{1}{\calD} & 0 \etbt.$

We seek $Z$ such that
$$ZAZ^{-1}=A,\quad ZC_1Z^{-1}= C_1,$$
and
$$Z\btbt 0 & \calD\\ -\frac{1}{\calD} & 0 \etbt Z^{-1} = \btbt 0&1\\
-1 & 0 \etbt.$$

But $Z$ must be of the form
$$Z= \btbt x&0\\ 0 & \frac{1}{x} \etbt, $$
so that
$$Z \btbt 0 & \cal{D}\\ -\frac{1}{\calD} &0 \etbt Z^{-1} = \btbt 0 &
x^2 \calD\\ -\frac{1}{x^2\calD} & 0 \etbt,$$
which can't be made equal to $C_1$.

\vskip 12pt

\centerline{\bf REFERENCES}

\begin{description}

\bibo{Be} P. Berard, ``Transplantation et Isospectralit\'e,'' Math.\
Ann.\ 292 (1992), pp. 547-560. 

\bibo{BGP} R. Brooks, R. Gornet, and P. Perry, ``Isoscattering
Schottky Manifolds,'' GAFA 10 (2000), pp. 307-326.

\bibo{BJP} D. Borthwick, C. Judge, and P. Perry, ``Determinants of
Laplacians and Isopolar Metrics on Surfaces of Infinite Area,''
preprint. 

\bibo{BP} R. Brooks and P. Perry, ``Isophasal Scattering Manifolds in
Two Dimensions,'' Comm.\ Math.\ Phys. 223 (2001), pp. 465-474.

\bibo{BT} R. Brooks and R. Tse, ``Isospectral Surfaces of Small
Genus,'' Nagoya Math.\ J.\ 107 (1987), pp. 13-24.

\bibo{Di} L. E. Dickson, {\it Linear Groups,} 1901.

\bibo{Su} T. Sunada, ``Riemannian Coverings and Isospectral
Manifolds,'' Ann.\ Math.\ 121 (1985), pp. 169-186.

\bibo{Zel} S. Zelditch, ``Kuznecov Sum Formulae and Szeg\"o Limit
Formulas on Manifolds,'' Comm.\ P.\ D.\ E. 17 (1992), pp. 221-260.

\end{description}

\end{document}